\documentclass[11pt]{amsart}
\usepackage{amssymb, latexsym}
\theoremstyle{plain}
\newtheorem{theorem}{Theorem}

\newtheorem*{2'}{Theorem 2'}
\newtheorem*{3'}{Theorem 3'}

\theoremstyle{remark}

\newtheorem*{Remark 1}{Remark 1}
\newtheorem*{Remark 2}{Remark 2}
\newtheorem*{Remark 3}{Remark 3}
\newtheorem*{Remark 4}{Remark 4}

\numberwithin{equation}{section}

\begin{document}

\title [Strange domain of attraction to Dickman distributions]
{On the strange domain of attraction to generalized Dickman distributions for  sums of independent random variables}

\author{Ross G. Pinsky}

%\noindent  pinsky@math.technion.ac.il\ \ \ \ tel: 972-4-829-4083\ \ \  fax: 972-4-829-3388

\dedicatory{Dedicated to the memory of Mark A. Pinsky (1940-2016)}

\address{Department of Mathematics\\
Technion---Israel Institute of Technology\\
Haifa, 32000\\ Israel}
\email{ pinsky@math.technion.ac.il}

\urladdr{http://www.math.technion.ac.il/~pinsky/}

\date{}

\begin{abstract}
Let $\{B_k\}_{k=1}^\infty, \{X_k\}_{k=1}^\infty$
all be
 independent   random variables. Assume that $\{B_k\}_{k=1}^\infty$ are $\{0,1\}$-valued Bernoulli random variables satisfying
$B_k\stackrel{\text{dist}}{=}\text{Ber}(p_k)$, with $\sum_{k=1}^\infty p_k=\infty$,
and assume that  $\{X_k\}_{k=1}^\infty$  satisfy
  $X_k>0$ and  $\mu_k\equiv EX_k<\infty$.
Let
$M_n=\sum_{k=1}^np_k\mu_k$,  assume that $M_n\to\infty$
and define the normalized sum of independent random variables
$W_n=\frac1{M_n}\sum_{k=1}^nB_kX_k$.
We give a general condition under which $W_n\stackrel{\text{dist}}{\to}c$, for some $c\in[0,1]$, and a general condition under which
$W_n$ converges weakly to a distribution from a family of distributions that includes the generalized Dickman distributions GD$(\theta),\theta>0$. In particular,
we obtain the following result, which reveals a strange domain of attraction to generalized Dickman distributions. Assume that  $\lim_{k\to\infty}\frac{X_k}{\mu_k}\stackrel{\text{dist}}{=}1$. Let $J_\mu,J_p$ be  nonnegative integers, let $c_\mu,c_p>0$ and let

\noindent $
\mu_n\sim c_\mu n^{a_0}\prod_{j=1}^{J_\mu}(\log^{(j)}n)^{a_j}, \
p_n\sim c_p\big({n^{b_0}\prod_{j=1}^{J_p}(\log^{(j)}n)^{b_j}}\big)^{-1}, \ b_{J_p}\neq0.
%\end{aligned}
$
If
$$
\begin{aligned}
&i.\ J_p\le J_\mu;\\
&ii.\  b_j=1, \ 0\le j\le J_p;\\
&iii.\  a_j=0, \ 0\le j\le J_p-1,\ \text{and}\ \ a_{J_p}>0,
\end{aligned}
$$
then $ \lim_{n\to\infty}W_n\stackrel{\text{dist}}{=}\frac1{\theta}\text{GD}(\theta),\ \text{where}\ \theta=\frac{c_p}{a_{J_p}}.
$
Otherwise,  $\lim_{n\to\infty}W_n\stackrel{\text{dist}}{=}c$, for some $c\in[0,1]$.
We also give an application to the statistics of the number of inversions in certain random shuffling schemes.
\end{abstract}

\maketitle
\section{Introduction and Statement of Results}

The Dickman function $\rho_1$ is  the unique function,  continuous on $(0,\infty)$,
 and satisfying the differential-delay equation
$$
\begin{aligned}
&\rho_1(x)=0,\ x\le0;\\
&\rho_1(x)=1,\ x\in(0,1];\\
&x\rho_1'(x)+\rho_1(x-1)=0, \ x>1.
\end{aligned}
$$
This function has an interesting role in number theory and probability, which we  describe briefly at the end of this section.
With a little work, one can show  that the Laplace transform
of $\rho_1$ is given by $\int_0^\infty \rho_1(x)e^{-\lambda x}dx=\exp(\gamma+\int_0^1\frac{e^{-\lambda x}-1}xdx)$,
where $\gamma$ is  Euler's constant.
From this it follows that $\int_0^\infty \rho_1(x)dx=e^\gamma$, and consequently, that $e^{-\gamma}\rho_1$ is a probability density on $[0,\infty)$.
We will call this probability distribution the \it Dickman distribution.\rm\ We  denote its density by $p_1(x)=e^{-\gamma}\rho_1(x)$, and
we denote by $D_1$ a random variable distributed according to the Dickman distribution. Differentiating the Laplace transform
$E\exp(-\lambda D_1)=\exp(\int_0^1\frac{e^{-\lambda x}-1}xdx)$ of $D_1$ at $\lambda=0$ shows that $ED_1=1$.
The distribution decays very rapidly; indeed, it is not hard to show that $p_1(x)\le \frac{e^{-\gamma}}{\Gamma(x+1)},\ x\ge0$ \cite{MV}.

In fact,  for all $\theta>0$, $\exp(\theta\int_0^1\frac{e^{-\lambda x}-1}xdx)$ is  the Laplace transform of a probability distribution. (We will prove this directly; however, this fact  follows from the theory of infinitely  divisible distributions, and shows that the distribution
in  question is infinitely divisible.)
This distribution has density
$p_\theta=\frac{e^{-\theta\gamma}}{\Gamma(\theta)}\rho_\theta$, where  $\rho_\theta$ satisfies the differential-delay equation
\begin{equation}\label{rhotheta}
\begin{aligned}
&\rho_\theta(x)=0,\ x\le0;\\
&\rho_\theta(x)=x^{\theta-1},\ 0< x\le 1;\\
&x\rho_\theta'(x)+(1-\theta)\rho_\theta(x)+\theta\rho_\theta(x-1)=0,\ x>1.
\end{aligned}
\end{equation}
We will call such distributions \it generalized Dickman distributions\rm\ and denote them by GD$(\theta)$. We denote by $D_\theta$ a random variable with the GD$(\theta)$ distribution.
Differentiating its Laplace transform at $\lambda=0$ shows that $ED_\theta=\theta$.
These distribution decays very rapidly; indeed, it is not hard to show that $p_\theta(x)\le \frac{C_\theta}{\Gamma(x+1)}, x\ge1$, for an appropriate constant $C_\theta$.
A fundamental fact about these distributions is that
\begin{equation}\label{fund}
D_\theta\stackrel{\text{dist}}{=}U^{\frac1\theta}(D_\theta+1),
\end{equation}
where $U$ is  distributed according to the uniform distribution on $[0,1]$, and $U$ and $D_\theta$ on the right
hand side above are independent.  From
\eqref{fund} it is immediate that
$$
D_\theta\stackrel{\text{dist}}{=}U_1^\frac1\theta+(U_1U_2)^\frac1\theta+(U_1U_2U_3)^\frac1\theta+\cdots,
$$
where $\{U_n\}_{n=1}^\infty$ are IID random variables  distributed according to the uniform distribution on $[0,1]$.
It will follow from the proof of Theorem \ref{1} below
 that $\exp(\theta\int_0^1\frac{e^{-\lambda x}-1}xdx)$ is the Laplace transform of a probability distribution.
  In section \ref{background} we will prove that  a random variable
with such a distribution satisfies \eqref{fund}, and  that
if a random variable satisfies \eqref{fund}, then it  has a density of the form $c_\theta\rho_\theta$, where
$\rho_\theta$ satisfies \eqref{rhotheta}. Thus, this paper is self-contained with regard to  all the above noted facts,
with the exception of the rate of decay and the value $\frac{e^{-\theta\gamma}}{\Gamma(\theta)}$ of the normalizing constant $c_\theta$ in $p_\theta$.
For more on these distributions, including a derivation of the normalizing constant,
see, for example, \cite{ABT} and \cite{PW}.

In fact, the scope of this paper   leads us to consider a more general family of distributions
than the generalized Dickman distributions. Let
$\mathcal{X}\ge0$ be a
random variable satisfying $E\mathcal{X}\le 1$.
Then, as we shall see, for $\theta>0$, there exists a distribution whose Laplace transform is
$\exp\big(\theta\int_0^1\frac{E\exp(-\lambda x\mathcal{X})-1}xdx\big)$. We will denote this distribution
by $GD^{(\mathcal{X})}(\theta)$ and we denote a random variable
with this distribution by $D_\theta^{(\mathcal{X})}$. (When $\mathcal{X}\equiv1$, we revert to the previous notation for generalized
Dickman distributions.)
Differentiating the Laplace transform at $\lambda=0$ shows that $ED_\theta^{(\mathcal{X})}=\theta E\mathcal{X}$.

Mimicking the proof of \eqref{fund} that  we give in section \ref{background} shows that
\begin{equation}\label{fundagain}
D_\theta^{(\mathcal{X})}=U^{\frac1\theta}(D_\theta^{(\mathcal{X})}+\mathcal{X}),
\end{equation}
where $U$ is  distributed according to the uniform distribution on $[0,1]$, and $U$, $D_\theta^{(\mathcal{X})}$
and $\mathcal{X}$
 on the right
hand side above are independent.  From
\eqref{fundagain} it is immediate that
$$
D_\theta\stackrel{\text{dist}}{=}\mathcal{X}_1U_1^{\frac1\theta}+\mathcal{X}_2(U_1U_2)^{\frac1\theta}+\mathcal{X}_3(U_1U_2U_3)^{\frac1\theta}+\cdots,
$$
where $\{U_n\}_{n=1}^\infty$ and $\{\mathcal{X}_n\}_{n=1}^\infty$ are mutually independent sequences of IID random variables,
with $U_1$  distributed according to the uniform distribution on $[0,1]$ and $\mathcal{X}_1$ distributed according to the
distribution  of $\mathcal{X}$.

It is known that the generalized Dickman distribution $GD(\theta)$ arises as the limiting distribution of
$\frac1n\sum_{k=1}^nkY_k$, where the $\{Y_k\}_{k=1}^\infty$ are independent random variables with $Y_k$ distributed according to the Poisson distribution
with parameter $\frac \theta k$ \cite{ABT}.
It is also known that the Dickman distribution $GD(1)$ arises as the limiting distribution of
$\frac1n\sum_{k=1}^n kY_k$ as $n\to\infty$, where the $\{Y_k\}_{k=1}^\infty$ are independent Bernoulli random variables
satisfying $P(Y_k=1)=1-P(Y_k=0)=\frac1k$.
Such behavior is in distinct contrast to the law of large numbers behavior of
 a ``well-behaved'' sequence of independent random variables $\{Z_k\}_{k=1}^\infty$ with finite first moments;  namely, that
 $\frac1{M_n}\sum_{k=1}^nZ_k$  converges in distribution to 1 as $n\to\infty$, where $M_n=\sum_{k=1}^nEZ_k$.

The purpose of this paper is to understand when  the law of large numbers fails and a
distribution from the family $\text{GD}^{(\mathcal{X})}(\theta)$  arises in its stead.
From the above examples, we see that generalized Dickman distributions sometimes arise as limits of normalized  sums from a  sequence $\{V_k\}_{k=1}^\infty$ of independent random variables which are
are non-negative and satisfy the following  three conditions: (i) $\lim_{k\to\infty}P(V_k=0)=1$,
(ii) $\lim_{k\to\infty}\frac{V_k|V_k>0}{E(V_k|V_k>0)}\stackrel{\text{dist}}{=}1$ and (iii) $\sum_{k=1}^\infty EV_k=\infty$.
(In the above examples, $kY_k$ plays the role of $V_k$.)
It turns out that these three conditions are very far from sufficient for a generalized Dickman distribution to arise. In fact, as we shall see in Theorem \ref{2} below,  such distributions  arise only
in  a  strange sequence of  very narrow windows of opportunity.

In light of the above discussion, we   will consider the following setting.
Let $\{B_k\}_{k=1}^\infty, \{X_k\}_{k=1}^\infty$  be mutually independent sequences of independent   random variables. Assume that $\{B_k\}_{k=1}^\infty$ are Bernoulli random variables satisfying:
\begin{equation}\label{Ber}
P(B_k=1)=1-P(B_k=0)=p_k\in[0,1),
\end{equation}
and assume that  $\{X_k\}_{k=1}^\infty$  satisfy:
\begin{equation}\label{Xn}
  X_k>0,\ \  \ \mu_k\equiv EX_k<\infty.
\end{equation}
Let
\begin{equation}\label{Mn}
M_n=\sum_{k=1}^np_k\mu_k,
\end{equation}
and define
\begin{equation}\label{W}
W_n=\frac1{M_n}\sum_{k=1}^nB_kX_k.
\end{equation}
We will be interested in the limiting behavior of $W_n$.
In order to avoid trivialities, we will assume that
\begin{equation}\label{nontriv}
\lim_{n\to\infty}M_n=\infty\ \ \text{and} \ \ \  \sum_{k=1}^\infty p_k=\infty,
\end{equation}
since otherwise $\sum_{n=1}^\infty B_kX_k$ is almost surely finite.

Note that for the  example brought with the $\text{Pois}(\frac\theta k)$-distribution, we have $p_k=1-e^{-\frac\theta k}$, $X_k$ is distributed according to $kY_k|\{Y_k>0\}$, where $Y_k$ has the Pois$(\frac\theta k)$ distribution,
$\mu_k=\frac\theta{1-e^{-\frac\theta k}}$
and $M_n=n\theta$. And for the example with the  $\text{Ber}(\frac1k)$-distribution, we have $p_k=\frac1k$, $X_k=k$ deterministically, $\mu_k=k$ and $M_n=n$.
In the first of these two examples, $\frac{X_k}{\mu_k}\stackrel{\text{dist}}{\to}1$,
and in the second one, $\frac{X_k}{\mu_k}\stackrel{\text{dist}}{=}1$ for all $k$.

Our first theorem  gives a  general condition for $W_n\stackrel{\text{dist}}{\to}c$ (which is the  law of large numbers if $c=1$), and
a general condition for convergence to a limiting distribution from the family of distributions
$\text{GD}^{(\mathcal{X})}(\theta)$.
Using this theorem, we can prove our second theorem, which reveals the strange domain of attraction to generalized Dickman distributions.
(Of course, we are using the term ``domain of attraction'' not in its classical sense, since our sequence of random variables, although independent,  are not identically distributed.)
Let $\delta_c$ denote the degenerate distribution at $c$.

\begin{theorem}\label{1}
Let $W_n$ be as in \eqref{W}, where $\{B_k\}_{k=1}^\infty$, $\{X_k\}_{k=1}^\infty$ and $M_n$ are as in
\eqref{Ber}-\eqref{Mn} and \eqref{nontriv}.

\noindent i. Assume that $\{\frac{X_k}{\mu_k}\}_{k=1}^\infty$ is uniformly integrable (which occurs automatically if $\lim_{k\to\infty}\frac{X_k}{\mu_k}\stackrel{\text{dist}}{=}1$).

\it a.
Assume also that
\begin{equation}\label{maxmueasy}
\lim_{n\to\infty}\frac{\max_{1\le k\le n}\thinspace\mu_k}{M_n}=0.
\end{equation}
Then
$$
\lim_{n\to\infty}W_n\stackrel{\text{dist}}{=}1.
$$

\it b.
 Assume  also that there exists a sequence $\{K_n\}_{n=1}^\infty$ such that
\begin{equation}\label{Kn}
\lim_{n\to\infty}\sum_{k=K_n+1}^n p_k=0,
\end{equation}
and
\begin{equation}\label{maxmu}
\lim_{n\to\infty}\frac{\max_{1\le k\le K_n}\thinspace\mu_k}{M_n}=0.
\end{equation}
If
\begin{equation}\label{cexists}
c\equiv\lim_{n\to\infty}\frac{M_{K_n}}{M_n}\ \text{exists},
\end{equation}
then
$$
\lim_{n\to\infty}W_n\stackrel{\text{dist}}{=}c.
$$
If \eqref{cexists} does not hold, then the distributions of $\{W_n\}_{n=1}^\infty$ form a tight sequence whose
set of accumulation points is $\{\delta_c: c\in A\}$, where
$A$ denotes the set of accumulation points of the sequence $\{\frac{M_{K_n}}{M_n}\}_{n=1}^\infty$.

\medskip

\noindent ii. Assume that there exists a random variable $\mathcal{X}$ such that
\begin{equation}\label{weakconvtoX}
\lim_{k\to\infty}\frac{X_k}{\mu_k}\stackrel{\text{dist}}{=}\mathcal{X}.
\end{equation}
Assume also that $\{\mu_k\}_{k=1}^\infty$ is increasing,  that $\lim_{k\to\infty}p_k=0$
and that there exist $\theta,L\in(0,\infty)$ such that
\begin{equation}\label{twocond}
\lim_{k\to\infty}\frac{p_k\mu_k}{\mu_{k+1}-\mu_k}=\theta,\ \ \ \lim_{k\to\infty}\frac{\mu_k}{M_k}=L.
\end{equation}
%Also assume that
%\begin{equation}\label{expmomentcond}
%\lim_{k_1\to\infty}\thinspace\lim_{n\to\infty}\thinspace\sup_{k_1\le k\le n}|E\exp(-\lambda\frac{X_k-\mu_k}{M_n})-1|=0.
%\end{equation}
Then
$$
\lim_{n\to\infty}W_n\stackrel{\text{dist}}{=}LD^{(\mathcal{X})}(\theta),
$$
where $D^{(\mathcal{X})}(\theta)$ is a random variable with the $\text{GD}^{(\mathcal{X})}(\theta)$ distribution.
\end{theorem}
\bf\noindent Remark 1.\rm\ In \eqref{twocond}, necessarily $L\le\frac1\theta$. Indeed, if $\{p_k\}_{k=1}^{\infty}$ and $\{\theta_k\}_{k=1}^\infty$ satisfy the conditions of part (ii),
and we choose $X_k=\mu_k$, then $W_n\stackrel{\text{dist}}{\to}LD_\theta$.
Since $EW_n=1$ and $ED_\theta=\theta$, it follows from Fatou's lemma
that $L\le \frac1\theta$.
In most  cases of interest, one has $L=\frac1\theta$.

\bf\noindent Remark 2.\rm\ By Fatou's lemma, the random variable $\mathcal{X}$ in part (ii) must satisfy
 $E\mathcal{X}\le 1$.

\noindent \bf Remark 3.\rm\  The   uniform integrability  of
$\{\frac{X_k}{\mu_k}\}_{k=1}^\infty$ in part (i)   occurs automatically
if $\lim_{k\to\infty}\frac{X_k}{\mu_k}\stackrel{\text{dist}}{=}1$, because
if a sequence  $\{Y_k\}_{k=1}^\infty$ of random variables satisfies $Y_k\stackrel{\text{dist}}{\to}Y$, and $E|Y_k|<\infty$,
then $E|Y_k|\to E|Y|$ is   equivalent to uniform integrability.

\noindent \bf Remark 4.\rm\ In the case that  $X_k=\mu_k$, or more generally, if  $EX_k\le C\mu_k^2$, for all $k$ and some $C>0$, then
$$
Var(W_n)\le\frac{C\sum_{k=1}^Np_k\mu_k^2}{M_n^2}=C\frac{\sum_{k=1}^np_k\mu_k^2}{(\sum_{k=1}^np_k\mu_k)^2}\le
C\frac{\sup_{1\le k\le n}\mu_k}{M_n}.
$$
Thus,  in this case part (i-a) follows directly from the second moment method.
 
\medskip

Using Theorem 1, we can prove the following theorem that exhibits the strange domain of attraction
to generalized Dickman distributions. Let $\log^{(j)}$ denote the $j$th iterate  of the logarithm, and
make the convention $\prod_{j=1}^0=1$.
\begin{theorem}\label{2}
Let $W_n$ be as in \eqref{W}, where $\{B_k\}_{k=1}^\infty$, $\{X_k\}_{k=1}^\infty$ and $M_n$ are as in
\eqref{Ber}-\eqref{Mn}. Assume also that  $\lim_{k\to\infty}\frac{X_k}{\mu_k}\stackrel{\text{dist}}{=}1$.
Let $J_\mu,J_p$ be  nonnegative integers, let $c_\mu,c_p>0$ and define
$$
\begin{aligned}
&\mu(x)= c_\mu x^{a_0}\prod_{j=1}^{J_\mu}(\log^{(j)}x)^{a_j}, \\
&p(x)= c_p\big({x^{b_0}\prod_{j=1}^{J_p}(\log^{(j)}x)^{b_j}}\big)^{-1},
\end{aligned}
$$
with $b_{J_p}\neq0$.
Assume that
$$
\begin{aligned}
&\mu_k\sim\mu(k),\ \ p_k\sim p(k);\\
&\mu_{k+1}-\mu_k\sim\mu'(k).
\end{aligned}
$$
Assume that the exponents $\{a_j\}_{j=0}^{J_\mu}, \{b_j\}_{j=0}^{J_p}$ have been chosen so that
\eqref{nontriv} holds.
If
\begin{equation}\label{3cond}
\begin{aligned}
&i.\ J_p\le J_\mu;\\
&ii.\  b_j=1, \ 0\le j\le J_p;\\
&iii.\  a_j=0, \ 0\le j\le J_p-1,\ \text{and}\ \ a_{J_p}>0,
\end{aligned}
\end{equation}
then
$$
\lim_{n\to\infty}W_n\stackrel{\text{dist}}{=}\frac1\theta D_\theta,\ \text{with}\ \theta=\frac{c_p}{a_{J_p}},
$$
where $D_\theta$ is a random variable with the GD$(\theta)$ distribution.

\noindent Otherwise,  $\lim_{n\to\infty}W_n\stackrel{\text{dist}}{=}c$, where $c\in\{0,1\}$.
To determine $c$, let
\begin{equation}\label{kappa}
\kappa_\mu=\min\{0\le j\le J_\mu:a_j\neq0\}\ \ \text{and}\  \ \kappa_p=\min\{0\le j\le J_p:b_j\neq1\}.
\end{equation}
If $\{0\le j\le J_\mu:a_j\neq0\}$ is not empty, $a_{\kappa_\mu}>0$ and either
$\{0\le j\le J_p:b_j\neq1\}$ is empty and $\kappa_\mu<J_p$, or
 $\{0\le j\le J_p:b_j\neq1\}$ is not empty and $\kappa_\mu<\kappa_p$, then
$c=0$; otherwise, $c=1$.
\end{theorem}
\bf\noindent Remark 1.\rm\ Note that if one chooses $\mu_k=\mu(k)$ and $p_k=p(k)$, then the condition
$\mu_{k+1}-\mu_k\sim \mu'(k)$ is always satisfied.

\noindent \bf Remark 2.\rm\ Theorem \ref{2} shows that to obtain a generalized Dickman distribution, $\{p_k\}_{k=1}^\infty$
in particular must be set in a very restricted fashion. For some intuition  regarding this phenomenon,
take the situation where $X_k=\mu_k$, and consider
 the sequence  $\{\sigma^2(W_n)\}_{n=1}^\infty$ of variances. This sequence
converges to 0 in the cases where $W_n$ converges to 1, converges to $\infty$ in the cases where $W_n$ converges to 0,
and converges to a positive number in  the cases where $W_n$ converges to a generalized Dickman distribution.

\medskip

We now state explicitly what  Theorem  \ref{2} yields in the cases $J_p=0,1$.
\medskip

\noindent $\bf J_p=0$.\rm\ We have
$$
p_n\sim \frac{c_p}{n^{b_0}},\  b_0>0, \ \ \mu_n\sim c_\mu n^{a_0}\prod_{j=1}^{J_\mu}(\log^{(j)}n)^{a_j}.
$$
In order that \eqref{nontriv} hold, we require
$b_0\le1$. We also require either: $a_0-b_0>-1$; or $a_0-b_0=-1$ and
 $a_1>-1$; or $a_0-b_0=a_1=-1$ and $a_2>-1$; etc.

\it\noindent  If $b_0=1$ and $a_0>0$, then
$$
\lim_{n\to\infty}W_n\stackrel{\text{dist}}{=} \frac1\theta D_\theta,\
\text{where}\ \theta=\frac{c_p}{a_0}.
$$
\noindent Otherwise,   $\lim_{n\to\infty}W_n\stackrel{\text{dist}}{=}1.$
\rm

%\noindent  Part (i) of the example  follows immediately from part (i) of the theorem with the choice $k_n=n$,
%and part (ii) of the example
% follows immediately from part (ii) of the theorem.

%\noindent \bf Remark.\rm\ Note that the constant $c_p$ in $p_k$ and the power $a$ in $\mu_k$ affect the limiting %distribution, but not the constant $c_\mu$ in $\mu_k$.
\medskip
%In light of the results in  Example I, we now fine tune around $p_n\sim\frac1n$.
%We begin by showing that  Example I-ii is robust with respect to lower order factors in $\mu_n$.

\medskip

\noindent $\bf J_p=1$.\rm\ We have
$$
p_n\sim\frac{c_p}{n^{b_0}(\log n)^{b_1}}, \ b_1\neq0,\  \ \ \mu_n\sim c_\mu n^{a_0}\prod_{j=1}^{J_\mu}(\log^{(j)}n)^{a_j}.
$$
In order that \eqref{nontriv} hold, we require
either $b_0=0$ and $b_1>0$, or $0<b_0<1$, or $b_0=1$ and $b_1\le1$. We also require
 either: $a_0-b_0>-1$; or $a_0-b_0=-1$ and
 $a_1-b_1>-1$; or $a_0-b_0=a_1-b_1=-1$ and $a_2>-1$; etc.

\it \noindent If $J_\mu\ge1$,  $b_0=b_1=1$, $a_0=0$ and $a_1>0$, then
$$
\lim_{n\to\infty}W_n\stackrel{\text{dist}}{=} \frac1\theta D_\theta,
\ \text{where}\ \theta=\frac{c_p}{a_1}.
$$
\medskip

\noindent If $b_0=1$ and $a_0>0$, then
$\lim_{n\to\infty}W_n\stackrel{\text{dist}}{=}0$.

\noindent Otherwise,
$\lim_{n\to\infty}W_n\stackrel{\text{dist}}{=}1$.
\medskip
\rm

\bf\noindent Remark.\rm\ In \cite{CS} and \cite{Pin}, where the GD$(1)$ distribution arises, one has
 $b_0=b_1=1, a_0=0, a_1=1,  c_p=c_\mu=1$.

The organization of the rest of the paper is as follows.
In section \ref{application} we
 use Theorems \ref{1} and  \ref{2} to investigate a  question raised in \cite{HT} concerning
 the statistics of the number of inversions in certain random shuffling schemes.
In sections \ref{pfthm1} and \ref{pfthm2} respectively we prove Theorems \ref{1} and \ref{2}.
Finally, in section \ref{background} we prove the basic facts about the Dickman distribution
and its density, as was promised earlier in this section.

As mentioned above,
we end this section with a little background concerning the Dickman function $\rho\equiv \rho_1$.
The Dickman function arises in probabilistic number theory in the context of so-called \it smooth \rm\ numbers; that is, numbers all of whose prime divisors are ``small.''
Let $\Psi(x,y)$ denote the number of positive integers less than or equal to $x$ with no
prime divisors greater than $y$. Numbers with no prime divisors greater than $y$ are  called $y$-\it smooth\rm\ numbers.
Then for $s\ge1$,
$\Psi(N,N^\frac1s)\sim N\rho(s)$, as $N\to\infty$.
%with the convergence uniform over bounded sets of $s$.
This result was first proved by Dickman in 1930 \cite{Dick}, whence the name of the function, with later refinements by de Bruijn \cite{deB1}.
See also \cite{MV} or \cite{Tene}.
Let $[n]=\{1,\ldots, n\}$ and
let $p^{+}(n)$ denote the largest prime divisor of $n$. Then Dickman's result states that
 the random variable $\frac{\log p^{+}(j)}{\log n}, j\in[n]$, on the probability space $[n]=\{1,\ldots, n\}$
  with the uniform distribution converges  in distribution as $n\to\infty$ to
the distribution whose distribution function is $\rho(\frac1x)$, $x\in[0,1]$, and whose density is $-\frac{\rho'(\frac1x)}{x^2}=\frac{\rho(\frac1x-1)}x,\  x\in[0,1]$.
It is easy to see that an equivalent statement of Dickman's
result is that the random variable $\frac{\log p^{+}(j)}{\log j}$, $j\in[n]$, on the probability space  $[n]$ with the uniform distribution converges in distribution as $n\to\infty$ to
the distribution whose distribution function is $\rho(\frac1x)$, $x\in[0,1]$,
We note that the length of the longest cycle of a uniformly random permutation of $[n]$, normalized by dividing by $n$,  also converges
to a limiting distribution whose distribution function is $\rho(\frac1x)$.
If instead of using the uniform measure on $S_n$, the set of permutations of $[n]$, one uses the Ewens sampling distribution on $S_n$, obtained
by giving each permutation $\sigma\in S_n$ the probability proportional to $\theta^{C(\sigma)}$, where $C(\sigma)$ denotes the number of cycles in $\sigma$, then
  the length of the longest cycle of such a  random permutation of $[n]$, normalized by dividing by $n$, converges to a limiting distribution whose distribution function
  is $\rho_\theta(\frac1x)$, $x\in[0,1]$.
  This distribution is also the distribution of the first coordinate of the Poisson-Dirichlet distribution PD$(\theta)$ (see \cite{ABT}).

  The examples in the above paragraph lead to limiting distributions where the Dickman function arises as a distribution function, not as a density as is the case with the GD$(\theta)$ distributions
  discussed in this paper. The GD$(\theta)$ distribution arises as a normalized limit in the context of certain natural probability measures that one can place on
  $\mathbb{N}$; see \cite{CS}, \cite{Pin}.

\section{An application to random permutations}\label{application}
We consider a setup that appeared in \cite{HT}, and which in
 the terminology of this paper can be described  as follows.
For each $k\in\mathbb{N}$, let $E_k\subset\{1,\ldots, k-1\}$.
Let $X_k$ be uniformly distributed on $E_k$, and let  $B_k\stackrel{\text{dist}}{=}\text{Ber}(\frac{|E_k|}{k})$.
So
$$
\mu_k=\frac1{|E_k|}\sum_{l\in E_k}l,\ \ \ p_k=\frac{|E_k|}k.
$$
Define
$$
I_n=\sum_{k=1}^n B_kX_k.
$$
We allow $E_k=\emptyset$, in which case $B_k=0$ and $X_k$ is not defined. In such a case, we define $B_kX_k=0$
and $\mu_k=0$. We always have $E_1=\emptyset$.

  Consider first the case that
 $E_k=\{1,\ldots, k-1\}$. Then $B_1X_1=0$ and for $2\le k\le n$,
 $B_kX_k$ is uniformly distributed over $\{0,1\ldots, k-1\}$. In this case, $I_n$ has the distribution of the number of inversions
in a uniformly random permutation
from $S_n$. (The authors in \cite{HT} have a typo and wrote $E_k=\{1,\ldots, k\}$  instead.)
To see this,  consider the following shuffling
procedure for  $n$ cards, numbered from 1 to $n$. The cards are to be
inserted in a row, one by one, in order of their numbers.
At step one,  card number 1 is set down. The number of inversions created by this step is zero, which is given by $B_1X_1$. At step $k$, for $k\in\{2,\ldots, n\}$,
card number $k$ is randomly inserted in the current row of cards, numbered 1 to $k-1$.
Thus, for any $j\in \{0,1,\ldots, k-1\}$, card number $k$ has probability $\frac1k$ of being placed
in the position with $j$ cards to its right (and $k-1-j$ cards to its left), in which  case this step
will have created $j$ new inversions, and this is represented by $B_kX_k$.
It is clear from the construction that the random variables $\{B_kX_k\}_{k=1}^n$ are independent.
Thus, $I_n$ indeed gives the number of inversions in a uniformly random permutation from $S_n$.
It  is well-known that the law of large numbers and the central limit theorem hold
for  $I_n$ in this case.

Consider now the general case that $E_k\subset\{1,\ldots, k-1\}$.
Then $I_n$ gives the number of inversions in a random permutation created
by a shuffling procedure in the same spirit as the above one. At step $k$, with probability $1-\frac{|E_k|}k$, card
number $k$
is inserted at the right end of the row, thereby creating no new inversions, and for each   $j\in E_k$,
with probability $\frac1k$ it is inserted
 in the position with $j$ cards to its right, thereby creating $j$ new inversions.

In particular, as a warmup consider the cases $E_k=\{1\}$ and   $E_k=\{k-1\}$,
$2\le k\le n$.
In each of these two cases, at step $k$, $2\le k\le n$, card number $k$ is inserted at the right  end of the row with
probability $1-\frac1k$.  In the first case, with probability $\frac1k$ card number $k$ is inserted immediately to the left
of the right most card, thereby creating one new inversion, while in the second case,
with probability $\frac1k$ card number $k$ is inserted at the left end of the row, thereby creating
$k-1$ new inversions.
In both  cases $\frac{X_n}{\mu_n}\stackrel{\text{dist}}{=}1$ for all $n$,  and
in both cases, $p_k=\frac1k$. In the first case, $\mu_k=1$ while in the second case, $\mu_k=k-1$.
Thus, in the first case, $M_n=\sum_{k=1}^n p_k\mu_k\sim\log n$, and in the second case,
$M_n\sim n$.
Therefore, it follows from Theorem \ref{1} or \ref{2} that in the first case $\frac{I_n}{\log n}$ converge in distribution
to 1, while  in the second case, $\frac{I_n}n$ converges in distribution  to GD$(1)$.

The authors of \cite{HT} ask which choices of $\{E_k\}_{k=1}^\infty$ lead to the  Dickman distribution and which
choices lead to the central limit theorem. Of course, the law of large numbers is a prerequisite
for the central limit theorem. The following theorem gives sufficient conditions for the law of large
numbers to  hold and sufficient conditions for convergence to a distribution from the
family  $\text{GD}^{(\mathcal{X})}(\theta)$.
In order to avoid trivialities, we need to assume that \eqref{nontriv} holds.
Recalling that $\mu_k=0$ when $|E_k|=0$, and that $\mu_k\ge1$ otherwise, note that
$$
M_n=EI_n=\sum_{k=1}^\infty\frac{|E_k|}k\mu_k\ge\sum_{k=1}^\infty\frac{|E_k|}k =\sum_{k=1}^\infty p_k.
$$
Thus, in the present context the requirement \eqref{nontriv} is
\begin{equation}\label{nontrivagain}
\sum_{k=1}^\infty\frac{|E_k|}k=\infty,
\end{equation}
which holds in particular if $E_k\neq\emptyset$ for all sufficiently large $k$.
\begin{theorem}\label{shufflethm}
Assume that \eqref{nontrivagain} holds.

\noindent i. Assume that at least one of the following conditions holds:

a. $\lim_{k\to\infty}|E_k|=\infty$ and $\{\frac{\mu_n}{\sum_{k=1}^n\frac{\mu_k}k}\}_{n=1}^\infty$ is bounded;

b. $\lim_{n\to\infty}\frac{\mu_n}{\sum_{k=1}^n\frac{\mu_k}k}=0$.

\noindent Then $\frac{I_n}{EI_n}\stackrel{\text{dist}}{\to}1$.

\noindent ii. Assume that $|E_k|=N\ge1$, for all large $k$,  and that $\frac{X_k}{\mu_k}\stackrel{\text{dist}}{\to} \mathcal{X}$.
 Also assume that  $\mu_k\sim\mu(k)$ and
$\mu_{k+1}-\mu_k\sim\mu'(k)$, where
$\mu(x)= c_\mu x^{a_0}\prod_{j=1}^{J_\mu}(\log^{(j)}x)^{a_j}$, with $a_0>0$.

\noindent Then $\frac{I_n}{EI_n}\stackrel{\text{dist}}{\to}\frac1\theta D_\theta^{(\mathcal{X})}$, with
$\theta=\frac{N}{a_0}$, where $D_\theta^{(\mathcal{X})}$ is a random variable with the $\text{GD}^{(\mathcal{X})}(\theta)$ distribution.

\end{theorem}

\bf\noindent Remark 1.\rm\ The condition on $\{\mu_k\}$ in  part (i-a) is just a very weak regularity requirement on its growth rate
(recall that $1\le \mu_k<k-1$).
The condition in part (i-b) is fulfilled if
$\mu_k\sim\mu(k)$ and
$j_{k+1}-j_k\sim\mu'(k)$, where
$\mu(x)= c_\mu \prod_{j=1}^{J_\mu}(\log^{(j)}x)^{a_j}$ with $J_\mu\ge0$.

\bf\noindent  Remark 2.\rm\ Note that  the random variable $\mathcal{X}$ in part (ii)  takes on no more than $N$ distinct values.

\begin{proof}
Assume first that the condition in part (i-a) holds.
We claim that since $\{\frac{\mu_n}{\sum_{k=1}^n\frac{\mu_k}k}\}_{n=1}^\infty$ is bounded, there exists a sequence of positive
integers
 $\{\gamma_n\}_{n=1}^\infty$ satisfying $\lim_{n\to\infty}\gamma_n=\infty$ and such that
 $\{\frac{\mu_n}{\sum_{k=\gamma_n+1}^n\frac{\mu_k}k}\}_{n=1}^\infty$ is also bounded.
 Indeed, assume to the  contrary. Then, in particular, $\{\mu_n\}_{n=1}^\infty$ is unbounded. Also,  since $\mu_k<k$, we have
 $\sum_{k=1}^n\frac{\mu_k}k<\gamma_n+\sum_{k=\gamma_n+1}^n\frac{\mu_k}k$, and it would follow
 that $\{\frac{\mu_n}{\gamma_n}\}_{n=1}^\infty$ is bounded for all sequences $\{\gamma_n\}_{n=1}^\infty$ satisfying
 $\lim_{n\to\infty}\gamma_n=\infty$, which is a contradiction.

 Let $\{\gamma_n\}_{n=1}^\infty$ be such a sequence. Then
$$
M_n=\sum_{k=1}^n\frac{|E_k|}k\mu_k\ge(\min_{k>\gamma_n}|E_k|)\sum_{k=\gamma_n}^n\frac{\mu_k}k.
$$
Thus, the condition in (i-a) guarantees that   \eqref{maxmueasy} holds.

Now assume that the condition in part (i-b) holds. Since
$M_n\ge\sum_{k=1}^n\frac{\mu_k}k$,  it follows again that \eqref{maxmueasy} holds.

 Thus, assuming  either (i-a) or (i-b), it follows  from
 part (i-a) of Theorem \ref{1} that $\frac{I_n}{EI_n}\stackrel{\text{dist}}{\to}1$.

Now assume that the condition in part  (ii) holds. Then  $p_k=\frac Nk$, for large $k$, and
$\mu_k\sim c_\mu k^{a_0}\prod_{j=1}^{J_\mu}(\log^{(j)}k)^{a_j}$,
with $a_0>0$. Thus,
$$
M_n=\sum_{k=1}^n\frac{|E_k|}k\mu_k\sim\frac {Nc_\mu}{a_0}n^{a_0}\prod_{j=1}^{J_\mu}(\log^{(j)}n)^{a_j},
$$
 and
 $\lim_{k\to\infty}\frac{\mu_k}{M_k}=\frac{a_0}N$.
Also, if the condition in part (ii) holds, then
\newline $\mu_{k+1}-\mu_k\sim a_0c_\mu k^{a_0-1} \prod_{j=1}^{J_\mu}(\log^{(j)}k)^{a_j}$.
Thus, $\lim_{k\to\infty}\frac{p_k\mu_k}{\mu_{k+1}-\mu_k}=\frac N{a_0}$.
We conclude from part (ii) of Theorem \ref{1} that $\frac{I_n}{EI_n}\stackrel{\text{dist}}{\to}\frac1\theta\text{GD}^{(\mathcal{X})}(\theta)$, where $\theta=\frac N{a_0}$.
\end{proof}

\section{Proof of Theorem \ref{1}}\label{pfthm1}
Since $EW_n=1$, for all $n$, the distributions of $\{W_n\}_{n=1}^\infty$ are tight.
Thus, since the random variables are nonnegative, it suffices to show  that  their
Laplace transforms $E\exp(-\lambda W_n)$ converge under the conditions
of part (i) to $\exp(-\lambda c)$, for the specified value of $c$,
and under the conditions of part (ii) to $\exp(\theta\int_0^1\frac{Ee^{-L\lambda x\mathcal{X}}-1}xdx)$,
which is the Laplace transform of $LD^{(\mathcal{X})}(\theta)$.
\medskip

\it \noindent Proof of part (i).\rm\
Note that part (i-a) is the particular case of part (i-b) in which one can choose $K_n=n$, and then
\eqref{cexists} holds with $c=1$. Thus, it suffices to consider part (i-b).
We have for $\lambda>0$,
\begin{equation}\label{LT}
\begin{aligned}
&E\exp(-\lambda W_n)=\\
&=\prod_{k=1}^nE\exp(-\frac\lambda{M_n} B_kX_k)=\prod_{k=1}^n\Big(1-p_k\big(1-E\exp(-\frac\lambda{M_n}X_k)\big)\Big)=\\
&\prod_{k=1}^{K_n}\Big(1-p_k\big(1-E\exp(-\frac\lambda{M_n}X_k)\big)\Big)\prod_{k=K_n+1}^n\Big(1-p_k\big(1-E\exp(-\frac\lambda{M_n}X_k)\big)\Big).
\end{aligned}
\end{equation}
Since
$$
\prod_{k=K_n+1}^n(1-p_k)\le\prod_{k=K_n+1}^n\Big(1-p_k\big(1-E\exp(-\frac\lambda{M_n}X_k)\big)\Big)\le1,
$$
it follows from
 assumption \eqref{Kn} that
\begin{equation}\label{Knassump}
\lim_{n\to\infty}\prod_{k=K_n+1}^n\Big(1-p_k\big(1-E\exp(-\frac\lambda{M_n}X_k)\big)\Big)=1.
\end{equation}

Applying the mean value theorem to $E\exp(-\frac\lambda{M_n}X_k)$ as a function of $\lambda$, and recalling that $\mu_k=EX_k$, we have
\begin{equation}\label{Taylor}
\frac{\lambda}{M_n}EX_k\exp(-\frac{\lambda}{M_n}X_k) \le1-E\exp(-\frac\lambda{M_n}X_k)\le\lambda\frac{\mu_k}{M_n}.
\end{equation}
In light of \eqref{maxmu} and the assumption that  $\{\frac{X_k}{\mu_k}\}_{k=1}^\infty$ is uniformly integrable, it follows that for all $\epsilon>0$, there exists an $n_\epsilon$ such that
\begin{equation}\label{fromwc}
\begin{aligned}
&\frac\lambda{M_n}EX_k\exp(-\frac{\lambda}{M_n}X_k)=\lambda\frac{\mu_k}{M_n} E\frac{X_k}{\mu_k}\exp(-\lambda\frac{\mu_k}{M_n}\frac{X_k}{\mu_k})\ge
(1-\epsilon)\lambda\frac{\mu_k}{M_n},\\
&  1\le k\le K_n,\ n\ge n_\epsilon.
\end{aligned}
\end{equation}
Thus, \eqref{Taylor} and \eqref{fromwc} yield
\begin{equation}\label{key}
(1-\epsilon)\lambda\frac{\mu_k}{M_n}  \le1-E\exp(-\frac\lambda{M_n}X_k)\le\lambda\frac{\mu_k}{M_n},\   1\le k\le K_n,\ n\ge n_\epsilon.
\end{equation}
Since for any $\epsilon>0$, there exists an $x_\epsilon>0$ such that $-(1+\epsilon)x\le \log (1-x)\le -x$, for $0<x<x_\epsilon$, it follows from \eqref{key} and  \eqref{maxmu} that
there exists an $n'_\epsilon$ such that
\begin{equation}\label{finallog}
\begin{aligned}
&-(1+\epsilon)\lambda p_k\frac{\mu_k}{M_n}\le\log\Big(1-p_k\big(1-E\exp(-\frac\lambda{M_n}X_k)\big)  \Big)\le-(1-\epsilon)\lambda p_k\frac{\mu_k}{M_n},\\
& 1\le k\le K_n, \   n\ge n'_\epsilon.
\end{aligned}
\end{equation}

From \eqref{finallog} we have
\begin{equation}\label{applylog}
\begin{aligned}
&-(1+\epsilon)\lambda\frac{\sum_{k=k_\epsilon}^{K_n}p_k\mu_k}{M_n}\le \log\prod_{k=1}^{K_n}\Big(1-p_k\big(1-E\exp(-\frac\lambda{M_n}X_k)\big)\Big)\le\\
&-(1-\epsilon)\lambda\frac{\sum_{k=k_\epsilon}^{K_n}p_k\mu_k}{M_n},\ n\ge n'_\epsilon.
\end{aligned}
\end{equation}
If
\begin{equation}\label{expectationcalc}
c\equiv\lim_{n\to\infty}\frac{M_{K_n}}{M_n}= \lim_{n\to\infty}\frac{\sum_{k=1}^{K_n}p_k\mu_k}{M_n}
\end{equation}
exists, then
 from \eqref{LT}, \eqref{Knassump}, \eqref{applylog} and \eqref{expectationcalc}, along with the fact that $\epsilon>0$ is arbitrary, we conclude that
$$
\lim_{n\to\infty}E\exp(-\lambda W_n)=\exp(-\lambda c),
$$
which proves that $\lim_{n\to\infty}W_n\stackrel{\text{dist}}{=}c$.
The rest of the results in part (i-b), concerning accumulation points, follow in the same manner.

\medskip

\noindent \it Proof of part (ii).\rm\
From \eqref{LT}, we have
\begin{equation}\label{logLT}
\log E\exp(-\lambda W_n)=\sum_{k=1}^n\log \Big(1-p_k\big(1-E\exp(-\frac\lambda{M_n}X_k)\big)\Big).
\end{equation}
Since by assumption $\lim_{k\to\infty}p_k=0$, for any $\epsilon>0$ there exists a $k_\epsilon$ such that
\begin{equation}\label{logapprox}
\begin{aligned}
&-(1+\epsilon)p_k\big(1-E\exp(-\frac\lambda{M_n}X_k)\big)\le \log \Big(1-p_k\big(1-E\exp(-\frac\lambda{M_n}X_k)\big)\Big)\le\\
&-p_k\big(1-E\exp(-\frac\lambda{M_n}X_k)\big), \ k\ge k_\epsilon.
\end{aligned}
\end{equation}

We now show that for any $\epsilon>0$ there exists a $k'_\epsilon$ such that
\begin{equation}\label{Xkmuk}
(1-\epsilon)E\exp(-\lambda\frac{\mu_k}{M_n}\mathcal{X})\le E\exp(-\frac\lambda{M_n}X_k)\le(1+\epsilon)E\exp(-\lambda\frac{\mu_k}{M_n}\mathcal{X}),\ k\ge k'_\epsilon.
\end{equation}
By assumption \eqref{twocond} and the  assumption that $\{\mu_n\}_{n=1}^\infty$ is increasing,
there exists a $C$ such that $\frac{\mu_k}{M_n}\le C$, for $1\le k\le n$ and $n\ge1$.
By assumption, $\frac{X_k}{\mu_k}\stackrel{\text{dist}}{\to}\mathcal{X}$. Without loss of generality, we assume that all of these random variables are defined on the same space and that
$\frac{X_k}{\mu_k}\to \mathcal{X}$ a.s.
For $\delta>0$, let
$$
A_{k;\delta}=\{\sup_{l\ge k}|\frac{X_l}{\mu_l}-\mathcal{X}|\le\delta\}.
$$
Then $A_{k;\delta}$ is increasing in $k$ and $\lim_{k\to\infty}P(A_{k;\delta})=1$.
We have
\begin{equation}\label{Adeltcompl}
\int_{A^c_{k;\delta}}\exp(-\frac\lambda{M_n}X_k)dP\le P(A^c_{k;\delta}),
\end{equation}
and
\begin{equation}\label{Adelt}
\begin{aligned}
&\exp(-\lambda C\delta)\int_{A_{k;\delta}}\exp(-\lambda\frac{\mu_k}{M_n}\mathcal{X})dP
\le\int_{A_{k;\delta}}\exp(-\lambda\frac{\mu_k}{M_n}\frac{X_k}{\mu_k})dP\le\\
& \exp(\lambda C\delta)\int_{A_{k;\delta}}\exp(-\lambda\frac{\mu_k}{M_n}\mathcal{X})dP.
\end{aligned}
\end{equation}
Now \eqref{Xkmuk}
follows from \eqref{Adeltcompl} and \eqref{Adelt}.

Letting $k^{''}_\epsilon=\max(k_\epsilon,k'_\epsilon)$, it follows from \eqref{logapprox} and \eqref{Xkmuk} that
\begin{equation}\label{goodlogapprox}
\begin{aligned}
&-(1+\epsilon)p_k\Big(1-(1-\epsilon)E\exp(-\lambda\frac{\mu_k}{M_n}\mathcal{X})\Big)\le \log \Big(1-p_k\big(1-E\exp(-\frac\lambda{M_n}X_k)\big)\Big)\le\\
&-p_k\Big(1-(1+\epsilon)E\exp(-\lambda\frac{\mu_k}{M_n}\mathcal{X})\Big),\ k\ge k^{''}_\epsilon.
\end{aligned}
\end{equation}
From \eqref{logLT} and \eqref{goodlogapprox} we have
\begin{equation}\label{forRiemsum}
\begin{aligned}
&-\sum_{k=k^{''}_\epsilon}^np_k(1+\epsilon)\Big(1-(1-\epsilon)E\exp(-\lambda\frac{\mu_k}{M_n}\mathcal{X})\Big)+o(1)\le\log E\exp(-\lambda W_n)\le\\
&-\sum_{k=k^{''}_\epsilon}^np_k\Big(1-(1+\epsilon)E\exp(-\lambda\frac{\mu_k}{M_n}\mathcal{X})\Big),\ \text{as}\ n\to\infty.
\end{aligned}
\end{equation}

Define $x_k^{(n)}=\frac{\mu_k}{M_n}$, $k^{''}_\epsilon\le k\le n$, and $\Delta_k^{(n)}=x_{k+1}^{(n)}-x_k^{(n)}=\frac{\mu_{k+1}-\mu_k}{M_n}$, $k^{''}_\epsilon\le k\le n-1$.
Then we have
\begin{equation}\label{Riemannsum}
\begin{aligned}
&\sum_{k=k^{''}_\epsilon}^np_k\Big(1-(1\pm\epsilon)E\exp(-\lambda\frac{\mu_k}{M_n}\mathcal{X})\Big)=\\
&\sum_{k=k^{''}_\epsilon}^n\frac{1-(1\pm\epsilon)E\exp(-\lambda x^{(n)}_k\mathcal{X})}{x_k^{(n)}}\Delta_k^{(n)}\big(p_k\frac{\mu_k}{\mu_{k+1}-\mu_k}\big).
\end{aligned}
\end{equation}
By assumption, $\{\mu_k\}_{k=1}^\infty$ is increasing; thus $\{x^{(n)}_k\}_{k=k^{''}_\epsilon}^n$ is a partition of $[\frac{\mu_{k^{''}_\epsilon}}{M_n},\frac{\mu_n}{M_n}]$.
By assumption, $\lim_{n\to\infty}\frac{\mu_{k^{''}_\epsilon}}{M_n}=0$ and $\lim_{n\to\infty}\frac{\mu_n}{M_n}=L$.
We now show
that the mesh,
$\max_{k^{''}_\epsilon\le k\le n-1}\Delta^{(n)}_k$,  of the partition
converges to 0 as $n\to\infty$.
Let $\Delta^{(n)}_{j_n}= \max_{k^{''}_\epsilon\le k\le n-1}\Delta^{(n)}_k$, where $k^{''}_\epsilon\le j_n\le n$.
Without loss of generality, assume either that $\{j_n\}$ is bounded or that $\lim_{n\to\infty}j_n=\infty$.
In the former case it is clear that $\max_{k^{''}_\epsilon\le k\le n-1}\Delta^{(n)}_k=\Delta^{(n)}_{j_n}=
\frac{\mu_{j_n+1}-\mu_{j_n}}{M_n}\stackrel{n\to\infty}{\to}0$.
Now consider the latter case.
From assumption \eqref{twocond} and the assumption that $\lim_{k\to\infty}p_k=0$, it follows that
$\lim_{n\to\infty}\frac{\mu_{n+1}-\mu_n}{M_n}=0$. Then we have
$$
\max_{k^{''}_\epsilon\le k\le n-1}\Delta^{(n)}_k=\Delta^{(n)}_{j_n}=
\frac{\mu_{j_n+1}-\mu_{j_n}}{M_n}=\frac{\mu_{j_n+1}-\mu_{j_n}}{M_{j_n}}\frac{M_{j_n}}{M_n}\le
\frac{\mu_{j_n+1}-\mu_{j_n}}{M_{j_n}}\stackrel{n\to\infty}{\to}0.
$$
Finally, we note that from \eqref{twocond} we have $\lim_{k\to\infty}p_k\frac{\mu_k}{\mu_{k+1}-\mu_k}=\theta$.
In light of these facts, along with \eqref{forRiemsum}, \eqref{Riemannsum} and the fact that $\epsilon>0$ is
arbitrary, it follows that
\begin{equation}
\lim_{n\to\infty}\log E\exp(-\lambda W_n)=\theta\int_0^L\frac{E\exp(-\lambda x\mathcal{X})-1}xdx=
\theta\int_0^1\frac{E\exp(-\lambda Lx\mathcal{X})-1}xdx,
\end{equation}

\hfill $\square$

\section{Proof of Theorem \ref{2}}\label{pfthm2}
We will assume that $J_p,J_\mu\ge1$ so that we can use a uniform notation, leaving it to the reader to verify
that the proof also goes through if $J_p$ or $J_\mu$ is equal to zero.

First assume that \eqref{3cond} holds. Then by the assumptions in the theorem,
\begin{equation*}
\begin{aligned}
&1\le J_p\le J_\mu;\\
&\mu_k\sim c_\mu\prod_{j=J_p}^{J_\mu}(\log^{(j)}k)^{a_j},\ a_{J_p}>0;\\
&p_k\sim c_p\big(x\prod_{j=1}^{J_p}\log^{(j)}k\big)^{-1};\\
&\mu_{k+1}-\mu_k\sim c_\mu a_{J_P}\thinspace\frac{(\log^{(J_P)}k)^{a_{J_p}-1}}{x\prod_{j=1}^{J_p-1}\log^{(j)}k}\thinspace
\prod_{j=J_p+1}^{J_\mu}(\log^{(j)}k)^{a_j}.
\end{aligned}
\end{equation*}
Thus,
\begin{equation*}
M_n=\sum_{k=1}^np_k\mu_k\sim c_\mu c_p\frac{(\log^{(J_p)}n)^{a_{J_p}}}{a_{J_p}}\prod_{j=J_p+1}^{J_\mu}(\log^{(  j)}n)^{a_j}.
\end{equation*}
Consequently,
\begin{equation}\label{secondcond}
\lim_{k\to\infty}\frac{\mu_k}{M_k}=\frac{a_{J_p}}{c_p}\ \ \text{and}\ \ \lim_{k\to\infty}\frac{p_k\mu_k}{\mu_{k+1}-\mu_k}=\frac{c_p}{a_{J_p}}.
\end{equation}
Thus, from part (ii) of Theorem \ref{1} it follows that
$\lim_{n\to\infty}W_n\stackrel{\text{dist}}{=}\frac1{\theta}D_\theta$, where $\theta=\frac{c_p}{a_{J_p}}$.

Now assume that \eqref{3cond} does not hold. We need to show that $\{K_n\}_{n=1}^\infty$ can be defined so that
\eqref{Kn} and \eqref{maxmu} hold, and so that \eqref{cexists} holds with $c\in\{0,1\}$. We also have to show when $c=0$ and when $c=1$.
Recall the definitions in \eqref{kappa}.
If $\{0\le j\le J_\mu:a_j\neq0\}$ is empty, or if it is not empty and $a_{\kappa_\mu}<0$,  then $\{\mu_k\}_{k=1}^\infty$ is bounded. Therefore,
\eqref{Kn} and \eqref{maxmu} hold with $K_n=n$ and it follows from  part (i-a) of Theorem \ref{1} that
$\lim_{n\to\infty}W_n\stackrel{\text{dist}}{=}1$. Thus, from now on we assume that
$\{0\le j\le J_\mu:a_j\neq0\}$ is not empty and that $a_{\kappa_\mu}>0$. In order to use uniform notation, we will
assume that $\kappa_\mu>0$, leaving the reader to verify that the proof goes through if $\kappa_\mu=0$.
Thus, we have
\begin{equation}\label{muformprelim}
\mu_k\sim\prod_{j=\kappa_\mu}^{J_\mu}(\log^{(j)}k)^{a_j},\ \ \kappa_\mu\ge1,\ a_{\kappa_\mu}>0.
\end{equation}
In order to simplify notation, for the rest of this proof, we will let
$\mathcal{L}_l(k)$  denote a positive constant multiplied by a product of powers (possibly of varying sign) of iterated logarithms
$\log^{(j)}k$, where the smallest $j$  is strictly larger than $l$. The exact from of this expression may vary from line to line.
Sometimes we will need to distinguish between two such expressions in the same formula, in which case we will use the notation
$\mathcal{L}^{(1)}_l(k), \mathcal{L}^{(2)}_l(k)$.
Thus, we rewrite \eqref{muformprelim} as
\begin{equation}\label{muform}
\mu_k\sim(\log^{(\kappa_\mu)}k)^{a_{\kappa_\mu}   }\mathcal{L}_{\kappa_\mu}(k),\ \ \kappa_\mu\ge1,\ a_{\kappa_\mu}>0.
\end{equation}

If $\{0\le j\le J_p:b_j\neq1\}$ is empty,  then the second condition in \eqref{3cond} is fulfilled and we have
\begin{equation}\label{pkempty}
p_k\sim c_p\big(x\prod_{j=1}^{J_p}\log^{(j)}k\big)^{-1}.
\end{equation}
Since
we are assuming that \eqref{3cond} does not hold, at least one of the other
two conditions in \eqref{3cond} must fail.
This forces $\kappa_\mu\neq J_p$.
(Recall that we are assuming that
$\{0\le j\le J_\mu:a_j\neq0\}$ is not empty and that $a_{\kappa_\mu}>0$.)

Consider first the case that $\kappa_\mu>J_p$.
Then from \eqref{muform} and \eqref{pkempty} we have
\begin{equation}\label{Mnfirst}
\begin{aligned}
&M_n=\sum_{k=1}^np_k\mu_k\sim
&(\log^{(J_p+1)}n)(\log^{(\kappa_\mu)}n)^{a_{\kappa_\mu}}\mathcal{L}_{\kappa_\mu}(n),
\ \text{where}\ \kappa_\mu\ge J_p+1.
\end{aligned}
\end{equation}
From \eqref{muform} and  \eqref{Mnfirst} it follows that \eqref{Kn} and \eqref{maxmu} hold by choosing $K_n=n$. Thus, from part (i-a) of Theorem \ref{1}, $\lim_{n\to\infty}W_n\stackrel{\text{dist}}{=}1$.

Now consider the case $\kappa_\mu<J_p$.
Then from \eqref{muform} and \eqref{pkempty} we have
\begin{equation}\label{Mncase1}
\begin{aligned}
&M_n=\sum_{k=1}^np_k\mu_k\sim
% \sum_{k=1}^n
%\frac{\prod_{j=\kappa_\mu}^{J_\mu}(\log^{(j)}k)^{a_j}}
%{k\prod_{j=1}^{J_p}\log^{(j)}k}\sim\\
&(\log^{(\kappa_\mu)}n)^{a_{\kappa_\mu}}\mathcal{L}_{\kappa_\mu}(n),
\ \text{where}\ \kappa_\mu\le J_p-1,
\end{aligned}
\end{equation}
and for any $K_n$ satisfying $K_n\to\infty$ and $K_n\le n$, we have
\begin{equation}\label{nKn}
\sum_{k=K_n}^np_k\sim c_p\big(\log^{(J_p+1)}n-\log^{(J_p+1)}K_n\big)=c_p\log\frac{\log^{(J_p)}n}{\log^{(J_p)}K_n}.
\end{equation}
From \eqref{muform} and \eqref{Mncase1} we have
\begin{equation}\label{muMn}
\begin{aligned}
&\frac{\mu_{K_n}}{M_n}
\sim\Big(\frac{\log^{(\kappa_\mu)}K_n}{\log^{(\kappa_\mu)}n}\Big)^{a_{\kappa_\mu}}\frac{\mathcal{L}^{(1)}_{\kappa_\mu}(K_n)}{\mathcal{L}^{(2)}_{\kappa_\mu}(n)},\ \kappa_\mu\le J_p-1,\ a_{\kappa_\mu}>0;\\
&\frac{M_{K_n}}{M_n}
\sim\Big(\frac{\log^{(\kappa_\mu)}K_n}{\log^{(\kappa_\mu)}n}\Big)^{a_{\kappa_\mu}}\frac{\mathcal{L}^{(1)}_{\kappa_\mu}(K_n)}{\mathcal{L}^{(2)}_{\kappa_\mu}(n)},\ \kappa_\mu\le J_p-1,\ a_{\kappa_\mu}>0;\\
\end{aligned}
\end{equation}
As we explain in some detail below, since $\kappa_\mu<J_p$, we can choose $\{K_n\}_{n=1}^\infty$ so that
\begin{equation}\label{Knbutn}
\lim_{n\to\infty}\frac{\log^{(J_p)}K_n}{\log^{(J_p)}n}=1\ \ \text{and}\ \
\lim_{n\to\infty}\Big(\frac{\log^{(\kappa_\mu)}K_n}{\log^{(\kappa_\mu)}n}\Big)^{a_{\kappa_\mu}}\frac{\mathcal{L}^{(1)}_{\kappa_\mu}(K_n)}{\mathcal{L}^{(2)}_{\kappa_\mu}(n)}=0.
\end{equation}
From \eqref{muform} and \eqref{nKn}-\eqref{Knbutn}, we conclude that $\{K_n\}$ can be defined so that
\eqref{Kn} and \eqref{maxmu} hold, and so that \eqref{cexists} holds with $c=0$. This proves that
$\lim_{n\to\infty}W_n\stackrel{\text{dist}}{=}0$.

To explain \eqref{Knbutn}, note that  $\frac{\mathcal{L}^{(1)}_{\kappa_\mu}(K_n)}{\mathcal{L}^{(2)}_{\kappa_\mu}(n)}\le(\log^{(\kappa_\mu+1)}n)^A$, for some
$A>0$ and all large $n$. (Recall that the powers of the iterated logarithms in $\mathcal{L}_{k_\mu}^{(2)}$ can be negative.)
%\begin{equation}\label{Knbutnagain}
%\lim_{n\to\infty}\frac{\log^{(J_p)}K_n}{\log^{(J_p)}n}=1\ \ \text{and}\ \
%\lim_{n\to\infty}\frac{\log^{(\kappa_\mu)}K_n}{(\log^{(\kappa_\mu)}n)^\frac12}=0.
%\end{equation}
Thus, in place of the second limit in \eqref{Knbutn}, it suffices to show that
 $\delta_n\equiv \Big(\frac{\log^{(\kappa_\mu)}K_n}{\log^{(\kappa_\mu)}n}\Big)^{a_{\kappa_\mu}}(\log^{(\kappa_\mu+1)}n)^A\stackrel{n\to\infty}{\to}0$.
 We have
 $$
 \log^{(\kappa_\mu)}K_n=(\delta_n)^\frac1{a_{\kappa_\mu}}(\log^{(\kappa_\mu+1)}n)^{-\frac A{a_{\kappa_\mu}}}\log^{(\kappa_\mu)}n;
 $$
thus,
\begin{equation}\label{+1+2}
\frac{\log^{(\kappa_\mu+1)}K_n}{\log^{(\kappa_\mu+1)}n}=\frac{\log\delta_n}{a_{\kappa_\mu}\log^{(\kappa_\mu+1)}n}-\frac{A\log^{(\kappa_\mu+2)}n}{a_{\kappa_\mu}\log^{(\kappa_\mu+1)}n}+1.
\end{equation}
Defining $K_n$ by choosing $\delta_n=(\log^{(\kappa_\mu+1)}n)^{-1}$,
it follows from \eqref{+1+2} and the fact that $J_p\ge \kappa_\mu+1$ that
the two equalities in \eqref{Knbutn} hold.

We now consider the case that
 $\{0\le j\le J_p:b_j\neq1\}$ is not empty. Then in order to fulfill the second condition  in
\eqref{nontriv}, we have $b_{\kappa_p}<1$. We write
\begin{equation}\label{pknotempty}
p_k\sim c_p\big(x\prod_{j=1}^{\kappa_p-1}\log^{(j)}k\big)^{-1}\big(\log^{(\kappa_p)}k\big)^{-b_{\kappa_p}}\big(\prod_{j=\kappa_p+1}^{J_p}\log^{(j)}k\big)^{-b_j}.
\end{equation}
From  \eqref{muform} and \eqref{pknotempty} it follows that
$M_n=\sum_{k=1}^np_k\mu_k$
satisfies
\begin{equation}\label{Mncases}
M_n\sim\begin{cases}(\log^{(\kappa_\mu)}n)^{a_{\kappa_\mu}}\thinspace\mathcal{L}_{\kappa\mu}(n),\ \kappa_\mu<\kappa_p;\\
(\log^{(\kappa_p)}n)^{a_{\kappa_p}-b_{\kappa_p}+1}\thinspace\mathcal{L}_{\kappa_p}(n),\ \kappa_\mu=\kappa_p;\\
(\log^{(\kappa_p)}n)^{1-b_{\kappa_p}}\thinspace\mathcal{L}_{\kappa_p}(n),\ \kappa_\mu>\kappa_p,\end{cases}
\end{equation}
and from \eqref{pknotempty} it follows that
 for any $K_n$ satisfying $K_n\to\infty$ and $K_n\le n$,
\begin{equation}\label{nKnagain}
\begin{aligned}
&\sum_{k=K_n}^np_k\sim\\
&\frac{c_p}{1-b_{\kappa_p}}\Big[\big(\log^{(\kappa_p)}n\big)^{1-b_{\kappa_p}}\big(\prod_{j=\kappa_p+1}^{J_p}\log^{(j)}n\big)^{-b_j}
-\big(\log^{(\kappa_p)}K_n\big)^{1-b_{\kappa_p}}\big(\prod_{j=\kappa_p+1}^{J_p}\log^{(j)}K_n\big)^{-b_j}\Big].
\end{aligned}
\end{equation}
From \eqref{muform} and \eqref{Mncases} we have
\begin{equation}\label{muMnagain}
\frac{\mu_{K_n}}{M_n}\sim\begin{cases}\big(\frac{\log^{(\kappa_\mu)}K_n}{\log^{(\kappa_\mu)}n}\big)^{a_{\kappa_\mu}}\thinspace\frac{\mathcal{L}^{(1)}_{\kappa_\mu}(K_n)}{\mathcal{L}^{(2)}_{\kappa_\mu}(n)},\ \kappa_\mu<\kappa_p;\\
 \frac{\big(\log^{(\kappa_p)}K_n)^{a_{\kappa_p}}}{(\log^{(\kappa_p)}n\big)^{a_{\kappa_p}-b_{\kappa_p}+1}}
\thinspace \frac{\mathcal{L}^{(1)}_{\kappa_p}(K_n)}{\mathcal{L}^{(2)}_{\kappa_p}(n)},\ \kappa_\mu=\kappa_p;\\
\frac{\big(\log^{(\kappa_\mu)}K_n\big)^{a_{\kappa_\mu}}}{\big(\log^{(\kappa_p)}n\big)^{1-b_{\kappa_p}}}\thinspace\frac{\mathcal{L}^{(1)}_{\kappa_\mu}(K_n)}{\mathcal{L}^{(2)}_{\kappa_p}(n)},\ \kappa_\mu>\kappa_p.
 \end{cases}
\end{equation}

It is immediate \eqref{muform} and \eqref{muMnagain} that if $\kappa_\mu\ge \kappa_p$, then \eqref{Kn} and \eqref{maxmu} hold by choosing $K_n=n$. (For the case $\kappa_\mu=\kappa_p$, recall that $b_{\kappa_p}\in(0,1)$.)
Thus, from part (i-a) of Theorem \ref{1}, $\lim_{n\to\infty}W_n\stackrel{\text{dist}}{=}1$.

Now consider the case $\kappa_\mu< \kappa_p$. For simplicity, we will
assume that the higher order iterated logarithmic terms do not appear; that is, we will assume from
\eqref{Mncases}-\eqref{muMnagain}  that
\begin{equation}\label{easier}
\begin{aligned}
&\sum_{k=K_n}^np_k\sim\frac{c_p}{1-b_{\kappa_p}}\Big[\big(\log^{(\kappa_p)}n\big)^{1-b_{\kappa_p}}
-\big(\log^{(\kappa_p)}K_n\big)^{1-b_{\kappa_p}}\Big];\\
&\frac{\mu_{K_n}}{M_n}\sim\big(\frac{\log^{(\kappa_\mu)}K_n}{\log^{(\kappa_\mu)}n}\big)^{a_{\kappa_\mu}};\\
&\frac{M_{K_n}}{M_n}\sim\big(\frac{\log^{(\kappa_\mu)}K_n}{\log^{(\kappa_\mu)}n}\big)^{a_{\kappa_\mu}}.
\end{aligned}
\end{equation}
The additional logarithmic terms can be dealt with similarly to the way they were dealt with for
\eqref{Knbutn}, as explained in the paragraph following \eqref{Knbutn}.
Applying the mean value theorem to the function $x^{1-b_{\kappa_p}}$, we obtain
\begin{equation}\label{MVT}
\big(\log^{(\kappa_p)}n\big)^{1-b_{\kappa_p}}
-\big(\log^{(\kappa_p)}K_n\big)^{1-b_{\kappa_p}}
=\frac{(1-b_{\kappa_p})\log^{(\kappa_p)}\frac n{K_n}}{(\log^{(\kappa_p)}n^*)^{b_{\kappa_p}}},
\end{equation}
where $n^*\in(K_n,n)$.
Since $\kappa_\mu<\kappa_p$, we can choose $K_n\to\infty$ such that
$\lim_{n\to\infty}\frac{\log^{(\kappa_\mu)}K_n}{\log^{(\kappa_\mu)}n}=0$,
but $\lim_{n\to\infty}\log^{(\kappa_p)}\frac {K_n}n=1$. For such a choice of $\{K_n\}$, it follows
from \eqref{muform}, \eqref{easier} and  \eqref{MVT} that \eqref{Kn} and \eqref{maxmu} hold, and that \eqref{cexists} holds with $c=0$;
thus, $\lim_{n\to\infty}W_n\stackrel{\text{dist}}{=}0$.
\hfill$\square$

\section{Basic Facts Concerning Generalized Dickman Distributions}\label{background}
We proved in Theorem \ref{1} that $\exp(\theta\int_0^1\frac{e^{-\lambda x}-1}xdx)$ is in fact the Laplace transform of a probability distribution, which we have denoted
by GD$(\theta)$.
In particular, if we let $X_k=\mu_k=k$ and $p_k=\frac\theta k$, in which case $M_n=\sum_{k=1}^np_k\mu_k=\theta n$, then it follows from Theorem \ref{1}
that
\begin{equation}\label{dickconv}
\hat W_n\equiv\theta W_n=\frac1n\sum_{k=1}^n kB_k\stackrel{\text{dist}}{\to}D_\theta,
\end{equation}
 where $D_\theta\stackrel{\text{dist}}{\sim} GD(\theta)$.

We now prove \eqref{fund}.
 Let
$$
J^+_n=\max\{k\le n: B_k\neq0\},
$$
with $\max\emptyset\equiv0$. We write
\begin{equation}\label{Whatn}
\hat W_n\equiv\frac1 n\sum_{n=1}^n kB_k=\frac{J^+_n-1}n\Big(\frac1{J^+_n-1}\sum_{k=1}^{J^+_n-1}kB_k\Big)+
\frac{J^+_n}n,
\end{equation}
where  the first of the two summands  on the right hand side above is interpreted as equal to 0 if
$J^+_n\le1$.
%Since $P(J^+_n\le k)=\prod_{m=k+1}^n(1-\frac\theta n)$,  we have $J^+_n\to\infty$ a.s. as $n\to\infty$.
We have
\begin{equation}\label{JNasym}
\begin{aligned}
P(\frac{J^+_n} n\le x)=\prod_{k=[xn+1]}^n(1-\frac\theta k)\sim x^\theta,\ x\in(0,1).
\end{aligned}
\end{equation}
Also, by the independence of $\{B_k\}_{k=1}^\infty$, we have
\begin{equation}\label{selfsim}
\frac1{J^+_n-1}\sum_{k=1}^{J^+_n-1}kB_k\mid   \{J^+_n=k_0\}\stackrel{\text{dist}}{=}\frac1{k_0-1}\sum_{k=1}^{k_0-1}kB_k=\hat W_{k_0-1},\ k_0\ge2.
\end{equation}
Letting $n\to\infty$ in \eqref{Whatn} and using \eqref{dickconv}, \eqref{JNasym} and \eqref{selfsim}, we conclude
that \eqref{fund} holds,
where $U$ is a distributed according to the uniform distribution on $[0,1]$, $D_\theta\stackrel{\text{dist}}{\sim}\text{GD}(\theta)$ and $U$ and $D_\theta$ on the right
hand side  are independent.

We now use \eqref{fund} to show that the GD$(\theta)$ distribution has a density function $p_\theta$ satisfying  $p_\theta=c_\theta\rho_\theta$, for some $c_\theta>0$, where
$\rho_\theta$ satisfies \eqref{rhotheta}. Let $F_\theta(x)=P(D_\theta\le x)$ denote the distribution function for the GD$(\theta)$ distribution.
Then from \eqref{fund} we have
\begin{equation}\label{Ftheta}
\begin{aligned}
&F_\theta(x)=P(D_\theta\le x)=P(U^\frac1\theta(D_\theta+1)\le x)=\int_0^1P(D_\theta+1\le xy^{-\frac1\theta})dy=\\
&\int_0^1F_\theta(xy^{-\frac1\theta}-1)dy.
\end{aligned}
\end{equation}
For $x>0$, making the change of variables, $v=xy^{-\frac1\theta}-1$, we can rewrite \eqref{Ftheta} as
\begin{equation}\label{Fthetaagain}
F_\theta(x)=\theta x^\theta\int_{x-1}^\infty F_{D_\theta}(v)(1+v)^{-1-\theta}dv,\ x>0.
\end{equation}
From \eqref{Fthetaagain} and the fact that  $F_\theta(x)=0$, for $x\le0$, it follows that $F_\theta$ is continuous on $\mathbb{R}$.
Also, since $F_\theta(x)=0$, for $x\le0$, we have
$$
\int_{x-1}^\infty F_{D_\theta}(v)(1+v)^{-1-\theta}dv=\int_0^\infty F_{D_\theta}(v)(1+v)^{-1-\theta}dv,\ x\le 1.
$$
Consequently,
it follows from \eqref{Fthetaagain} that
$F_\theta(x)=C_\theta x^\theta$, for $x\in[0,1]$, where $C_\theta=\theta\int_0^\infty F_{D_\theta}(v)(1+v)^{-1-\theta}dv$.
From this and \eqref{Fthetaagain} it follows that $F$ is differentiable on $(0,1)$ and on $(1,\infty)$,   and that, letting $p_\theta=F_\theta'$,
\begin{equation}\label{density01}
p_\theta=c_\theta x^{\theta-1},\ 0<x<1,\ c_\theta=\theta^2\int_0^\infty F_{D_\theta}(v)(1+v)^{-1-\theta}dv,
\end{equation}
and
\begin{equation}\label{density}
\begin{aligned}
&p_\theta(x)=\theta^2x^{\theta-1}\int_{x-1}^\infty F_{D_\theta}(v)(1+v)^{-1-\theta}dv-\theta x^{-1} F_\theta(x-1)=\\
&\frac\theta x (F_\theta(x)-F_\theta(x-1)), \ x>1.
\end{aligned}
\end{equation}
From \eqref{density}, it follows that $p_\theta$ is differentiable on $x>1$, and that
$(xp_\theta(x))'=\theta\big(p_\theta(x)-p_\theta(x-1)\big)$, for  $x>1$, or equivalently,
\begin{equation}\label{densityderiv}
xp'_\theta(x)+(1-\theta)p_\theta(x)+\theta p_\theta(x-1)=0,\ x>1.
\end{equation}
From \eqref{density01} and \eqref{densityderiv} we conclude that $p_\theta(x)=c_\theta\rho_\theta$, where $\rho_\theta$ satisfies \eqref{rhotheta}.
Integrating by parts in the formula for $c_\theta$ in \eqref{density01} shows that
$$
c_\theta=\theta\int_0^\infty(1+v)^{-\theta}p_\theta(v)dv=\theta E(1+D_\theta)^{-\theta}.
$$

\end{document}